\documentclass[10pt,a4paper,reqno]{amsart}

\usepackage{amsfonts,amsthm,latexsym,amsmath,amssymb,amscd,amsmath,epsf}

\makeatletter

\newtheorem{theorem}{Theorem}[section]
\newtheorem{theor}{Theorem}
\newtheorem*{coro}{Corollary}

\newtheorem{cor}[theorem]{Corollary}
\newtheorem{lemma}[theorem]{Lemma}
\newtheorem{ex}[theorem]{Example}
\newtheorem{st}[theorem]{Proposition}
\theoremstyle{definition}
\newtheorem{de}[theorem]{Definition}
\newtheorem{re}[theorem]{Remark}

\newcommand{\bC}{\mathbb C}
\newcommand{\bR}{\mathbb R}

\makeatother

\begin{document}

\title[Preserving positive polynomials and beyond]{Preserving positive polynomials and beyond}
\author[J.~Borcea]{Julius Borcea}
\address{Department of Mathematics, Stockholm University, SE-106 91
Stockholm, Sweden}
\email{julius@math.su.se}

\author[A.~Guterman]{Alexander Guterman}
       \address{Faculty of Algebra, Department of Mathematics and Mechanics, Moscow State University,
       119991, GSP-1, Moscow, Russia}
       \email{guterman@list.ru}

\author[B.~Shapiro]{Boris Shapiro}
\address{Department of Mathematics, Stockholm University, SE-106 91
Stockholm,
         Sweden}
\email{shapiro@math.su.se}



\begin{abstract}

Following the classical approach of P\'olya-Schur theory \cite {PS} we initiate  in this paper  the study  of   linear operators acting on $\bR[x]$ and preserving either the set of positive univariate polynomials or similar sets of non-negative and elliptic polynomials. 
\end{abstract}
\maketitle

\tableofcontents

\section{Introduction and main results}

Let ${\mathbb R}[x]$ denote the ring of univariate polynomials with real coefficients and denote by  ${\mathbb R}_n[x]$  its linear subspace  consisting of all polynomials of degree less than or equal to $n$.

In what follows we will discuss the following five important types of univariate polynomials:

\begin{de}
A polynomial $p(x)\in {\mathbb R}[x]$ is called 

{\em - hyperbolic\/}, if all its roots are real; 

 {\em - elliptic\/}, if it does not have reals roots; 

 {\em - positive\/}, if $p(x)>0$ for all $x\in {\mathbb R}$; 

 {\em - non-negative\/}, if $p(x)\ge 0$ for all $x\in {\mathbb R}$; 

 - a {\em sum of squares\/}, if there is a positive integer $k$ and there are polynomials $p_1(x),\ldots, p_k(x)\in {\mathbb R}[x]$ such that $p(x)=p^2_1(x)+\ldots + p^2_k(x)$. 
\end{de}

Note that the term ``elliptic''  is sometimes  used to define other types of polynomials, see, e.g., \cite{Lom,Pak}. The set of non-negative polynomials is classically  compared with the set of  sums of squares which is a subset of the latter. Moreover, a well-known result claims that in the univariate case these two classes coincide, see, e.g., \cite[p.~132]{W}. 

\begin{st}  \label{L:SumSq}
A polynomial $p(x)\in {\mathbb R}[x]$ is non-negative if and only if there exist $p_1(x),p_2(x)\in {\mathbb R}[x]$ such that $p(x)=p_1^2(x)+p_2^2(x)$.
\end{st}



\begin{re}
Note that the situation is quite different for polynomials in several variables. In particular,  even in 2 variables not all non-negative polynomials  can be represented as sums of squares. One of the simplest examples of this kind is the polynomial $p(x,y)=x^2y^2(x^2+y^2-3)+1$ which is non-negative but can not be represented as the sum of squares, see~\cite{Mo} for details. In general, this topic is related to the  Hilbert 17-th problem, see~\cite{PD}.
\end{re}

\begin{de}
Let $V$ denote either ${\mathbb R}_n[x]$ or ${\mathbb R}[x]$. We say that a map $\Phi:V \to V$ {\em preserves\/} a certain set $M\subset V$ if for any $p(x)\in M$ its image $\Phi(p(x))$ belongs to $M$. 
\end{de}

In this paper we study linear operators on ${\mathbb R}[x]$ or ${\mathbb R}_n[x]$ which preserve one of the classes of polynomials introduced above. Namely, we  call a linear operator acting on 
$\bR[x]$ or ${\mathbb R}_n[x]$ a {\em hyperbolicity-, ellipticity-, positivity-, non-negativity-preserver}  if it preserves the sets of hyperbolic, elliptic, positive, non-negative polynomials respectively. The classical case of (linear) hyperbolicity-preservers which are diagonal in the monomial basis of $\bR[x]$ was thoroughly studied  about a century ago by P\'olya and Schur \cite{PS}. Its substantial generalizations both in the univariate and the multivariate cases can be found in \cite{BB,BBS1,BBS2}.

Following the set-up of \cite {PS}  we concentrate below on the remaining three classes  of preservers (restricting our attention mainly to linear ordinary differential operators of finite order, see Remark~\ref{rmk}). In short, it turns out that there are much fewer such linear operators  than those preserving hyperbolicity. 
More precisely, our  two main results are as follows.

\begin{theor} \label{TT1}
Let $U_Q:{\mathbb R}[x]\to {\mathbb R}[x]$ be a linear ordinary differential operator of order $k\ge 1$ with polynomial coefficients $Q=(q_0(x), q_1(x),\ldots, q_k(x))$, $q_i(x)\in {\mathbb R}[x]$, $i=0,\ldots,k$, $q_k(x)\not\equiv 0$, i.e., 
\begin{equation}\label{BasOper}
 U_Q=q_0(x)+q_1(x) \frac{d}{dx}+q_2(x)\frac{d^2}{dx^2}+\ldots + q_k(x) \frac{d^k}{dx^k}.
 \end{equation}
Then for any coefficient sequence $Q$ 
the operator $U_Q$ does not preserve the set of non-negative
(resp., positive or elliptic)
 polynomials of degree $2k$.
\end{theor}

\begin{coro}
There are no   linear ordinary differential operators of positive finite order which preserve the set of 
non-negative 
(resp., positive  or elliptic)
polynomials in ${\mathbb R}[x]$.
\end{coro}

\begin{re} Notice that by contrast with the above situation there are many hyperbolicity-preservers which are finite order linear differential operators with polynomial coefficients. In fact, such examples exist even among operators with constant coefficients, see Remark~\ref{rmk} and \cite {BBS1}. 

\end{re}

\begin{re}\label{rmk}  Any linear operator on $\bR[x]$ and $\bC[x]$ can be represented as  a linear ordinary differential operator of, in general, {\em infinite} order, i.e., as a formal power series in $\frac{d}{dx}$ with polynomial coefficients, see, e.g., \cite{BBS1}.  Thus the subclass of finite order linear differential operators, i.e., those belonging to the Weyl algebra $\mathcal A_1$ is a natural object of study. Note that unlike the case of finite order operators  there exist plenty of linear differential operators of infinite order which preserve positivity. Apparently, the simplest example of this kind is 
\begin{equation}\label{exx}
\left(1-\frac {d}{dx}\right)^{-1}=1+\frac{d}{dx}+\frac{d^2}{dx^2}+\ldots
\end{equation}
More generally, the inverse of any finite order differential operator with constant coefficients and positive constant term  whose symbol is a hyperbolic polynomial yields an example of such an operator. \end{re}

Luckily 
the case of infinite order linear differential operators with constant coefficients can be handled completely. Namely, slightly generalizing a one hundred years old result of Remak \cite{Re} and  Hurwitz \cite{Hu} (see also Problem 38 in \cite[Ch.~7]{PSz})  one obtains the following statement.

\begin{theor} \label{TC:Fd}
Let  $\alpha=(\alpha_0,\alpha_1,\ldots,\alpha_k, ..., )$ be an infinite  sequence of real numbers.  Consider  the  infinite order linear ordinary differential operator   
\begin{equation}\label{ConsOper}
U_{\alpha}=\alpha_0+\alpha_1\frac{d}{dx}+\alpha_2\frac{d^2}{dx^2}+\ldots + \alpha_k \frac{d^k}{dx^k}+ ...
\end{equation}
with constant coefficients. Then the operator $U_{\alpha}$ 
 preserves positivity (resp., non-negativity) if and only if one of the following two equivalent conditions holds:
 \begin{enumerate}
 \item   for any positive (resp., non-negative) polynomial $p(x)=a_{k}x^{k}+\ldots + a_1x+a_0$ one has  that 
 $$ U_\alpha(p(0))=a_0\alpha_0+a_1\alpha_1+\ldots + k! a_{k} \alpha_{k}>0 \ (\mbox{resp., }\ge 0); 
 $$
 \item 	the following infinite  Hankel matrix 
 $$ \begin{pmatrix}  \alpha_0 & 1!\alpha_1 & 2!\alpha_2 & \ldots & l!\alpha_l &\ldots \\ 1!\alpha_1 & 2!\alpha_2 & 3!\alpha_3 & \ldots & (l+1)!\alpha_{l+1}&\ldots \\  2!\alpha_2 & 3!\alpha_3& 4!\alpha_4 & \ldots & (l+2)!\alpha_{l+2}&\ldots  \\ \vdots & \vdots & \vdots & \ddots& \vdots &\vdots \\ l!\alpha_l & (l+1)!\alpha_{l+1} & (l+2)!\alpha_{l+2} & \ldots & (2l)!\alpha_{2l} &\ldots \\
 \ldots & \ldots & \ldots & \ldots & \ldots & \ldots \end{pmatrix}$$ 
is positive definite (resp., positive semi-definite), i.e. all its principal minors  are positive (resp. non-negative).    
 \end{enumerate}
 \end{theor}

\medskip

To illustrate the latter result notice that for the operator~(\ref{exx}) above one has $1=\alpha_0=\alpha_1=\alpha_2=\alpha_3=...$ and $\Delta_l=\Pi_{j=1}^l (j!)^2,\; l=0,1,...$ where $\Delta_l$ is the corresponding $(l+1)\times (l+1)$ principal minor, see \cite {Ehr}.  

  \begin{re}
The major remaining challenge in this area is to classify {\em all} positivity-preservers.  We finish our introduction with this question. 
\end{re} 

\medskip 
\noindent
{\bf Problem 1.} Find a complete classification of positivity-preservers. 
\medskip 
 
We also state a more concrete and (hopefully) simpler question. 

\medskip 
\noindent  
 {\bf Problem 2.} 
Is it true that any positivity-preserver which is an infinite order linear differential operator with constant coefficients has  a hyperbolicity-preserver as its inverse?

 \medskip 
 \noindent 
 {\bf Acknowledgments.} The authors  are grateful to Petter Br\"anden for important references. The second author is sincerely grateful to the Swedish Royal Academy of Sciences and the Mittag-Leffler Institute for supporting his visit to Stockholm in Spring 2007 when a substantial part of this project was carried out. The first and third authors would like to thank the  American Institute of Mathematics for its hospitality in May 2007.

\section{Some preliminaries on the considered classes of preservers }

Below we discuss the relationships between the classes of ellipticity-, positivity-, and  non-negativity-%
preservers. 
As we mentioned in the introduction the set of all univariate non-negative polynomials coincides with the set of   sums of squares and therefore  linear preservers of the latter set do not require separate consideration. On the other hand, it is obvious  that the sets of elliptic, positive, and non-negative polynomials are distinct.  In this section we answer the  question about how different are the corresponding sets of ellipticity-, positivity- and non-negativity-preservers, respectively, see Theorems~\ref{Lnn1} and \ref{Lnn} below. 

We start with the following lemma showing  that the assumption that  a linear operator  $\Phi$ is a  non-negativity-preserver   is quite strong. 
\begin{lemma} \label{L:Phi}
Let $\Phi: {\mathbb R}[x]\to {\mathbb R}[x]$  be a linear operator preserving the set of non-negative polynomials. If $\Phi(1)\equiv  0$ then  $\Phi\equiv 0$.
\end{lemma}

\begin{proof}
Assume that $\Phi(1)\equiv 0$. First we show that  for any polynomial $p(x)=a_n x^n+\ldots$ of even degree $n=2m$ one has  that if $a_n>0$ then $\Phi(p(x))$ is non-negative  and  if $a_n<0$ then $\Phi(p(x))$ is non-positive. Indeed, if $n$ is even and $a_n>0 $ then  $p(x)$ has a global minimum, say $M$. Thus $p(x)+|M|\ge 0$ for all $x\in {\mathbb R}$. Therefore, $\Phi(p(x)+|M|)\ge 0$ for all $x\in {\mathbb R}$. However, by  linearity and the assumption $\Phi(1)\equiv 0$ we get that $$\Phi(p(x))=\Phi(p(x))+|M|\Phi(1)=\Phi(p(x)+|M|)\ge0$$ for all $x\in {\mathbb R}$. For $a_n<0$ the result follows  by  linearity.

Now let us  show that $\Phi(1)\equiv 0$ implies  that $\Phi\equiv 0$. 
Let $q(x)$ be a polynomial such that $\Phi(q(x))=a_nx^n+\ldots +a_i x^i$ with  the smallest possible non-negative value of $i$ such that $a_i\ne 0$. Let $p(x)$ be a monic real polynomial of even degree satisfying the  condition  $\deg(q(x))<\deg(p(x))$. Thus $p(x)+\mu q(x)$ is monic for any $\mu\in {\mathbb R}$. The above argument shows that  $\Phi(p(x)+\mu q(x))$ is non-negative for all $\mu\in {\mathbb R}$. Notice that our choice  of $q(x)$ implies that the polynomial $\Phi(p(x))$ has vanishing coefficients at the degrees $0,\ldots,i-1$. Hence $\Phi(p(x))=b_l x^l + \ldots + b_ix^i$ for some positive integer $l$ and some coefficients $b_l,\ldots, b_i\in {\mathbb R}$. Then for any given $\mu$ there exists $g_\mu(x)\in {\mathbb R}[x]$ such that $\Phi(p(x))+\mu \Phi(q(x))=x^ig_\mu(x)$. Obviously,   the constant term of $g_\mu(x)$  equals  $b_i+\mu a_i$. Since $a_i\neq 0$ there exists $\mu_0\in {\mathbb R}$ such that $b_i+\mu_0a_i=g_{\mu_0}(0)<0$ and by continuity it follows that there exists a neighborhood $N(0)$ of the origin such that $g_{\mu_0}(x)<0$ for all $x\in N(0)$. Therefore, there exists $0\ne x_0\in N(0)$ such that $x_0^i>0$, hence $x_0^ig_{\mu_0}(x_0)<0$. This contradicts the assumption that $\Phi$ is a non-negativity-preserver. Thus $\Phi(q(x))$ has a vanishing term of degree $i$, which  contradicts the choice of $q(x)$. We deduce that $\Phi(q(x))\equiv 0$ for all $q(x)\in {\mathbb R}[x]$.
 \end{proof}

\begin{theorem} \label{Lnn2}
Let $\Phi: {\mathbb R}[x]\to {\mathbb R}[x]$  be a linear operator. Then the following conditions are equivalent:
\begin{enumerate}
\item $\Phi$ preserves the set of elliptic polynomials; 
\item either $\Phi$ or $-\Phi$ preserves the set of positive polynomials. 
\end{enumerate}
Also each of these conditions implies that 
\begin{enumerate}
\setcounter{enumi}{2}
\item either $\Phi$ or $-\Phi$ preserves the set of non-negative polynomials. 
\end{enumerate}
\end{theorem}

\begin{proof}
Note that the identically zero operator satisfies neither  condition (1) nor condition (2). Therefore,  we will assume that $\Phi\not\equiv 0$.

$(1) \Rightarrow (2)$. Assume that $\Phi$ preserves the set of elliptic polynomials and that  neither $\Phi$ nor $-\Phi$ preserves positivity. In other words, since $\Phi$ is an ellipticity-preserver this means that there exist positive  polynomials $p(x),q(x)\in {\mathbb R}[x]$  such that $\Phi(p(x))>0$ and $\Phi(q(x))<0$ for all $x\in \bR$. Note that no elliptic polynomials can be annihilated by  $\Phi$ since $0$ is not an elliptic polynomial.  
We consider the following two subcases:

\noindent
A. There exist two positive polynomials $p(x)$, $q(x)$ as above such that $\deg \Phi(p(x)) \ne \deg \Phi(q(x))$. Wlog we can assume that  $\deg \Phi(p(x)) > \deg \Phi(q(x))$. Since $\Phi(p(x))$ is a positive polynomial it has even degree and positive leading coefficient. Thus for any $\mu\in {\mathbb R}$ the polynomial  $\Phi(p(x))+\mu\Phi(q(x))$ has the same properties, i.e., is of even degree and has positive leading coefficient. Hence  there exists $x_0(\mu)\in {\mathbb R}$ such that $\Phi(p(x))+\mu\Phi(q(x))>0$ for all $x$ with $|x|>x_0(\mu)$.  

Now  set $y_0:=\Phi(p)(0)$ and  $z_0:=\Phi(q)(0)$. Obviously, $y_0>0$ and $z_0<0$ since $\Phi(p)$ is positive and  $\Phi(q)$ is negative. Let $\mu_0=\frac {2y_0}{-z_0}>0$. Then $p(x)+\mu_0 q(x)$ is positive since it is the sum of two positive polynomials. At the same time  for its image we have that $\Phi( p+\mu_0 q)(x)>0$ for $x>x_0(\mu_0)$. However, at the origin one has $$\Phi( p+\mu_0 q)(0)=y_0+\frac{2y_0}{-z_0}z_0=y_0-2y_0=-y_0<0,$$
so by continuity  $\Phi( p+\mu_0 q)(x)$ must have at least one real zero, which is a contradiction.

\noindent
B. It remains to consider the case when the images of all positive polynomials have the same degree, say $m$. Let $\Phi(p(x))=a_mx^m+\ldots,$ $\Phi(q(x))=b_mx^m+\ldots$. Since $\Phi(p(x))>0$ it follows that $a_m>0,$  and since $\Phi(q(x))<0$ one has $b_m<0$. Thus the polynomial $-b_mp(x)+a_mq(x)$ is positive. However, its image is of the  degree less than $m$, which is a contradiction.

 $(2) \Rightarrow (1)$. If $\Phi$ is a positivity-preserver  then by linearity $\Phi$ is also a negativity-preserver, and thus $\Phi$ preserves  the set of elliptic polynomials as well.

 $(2) \Rightarrow (3)$. Assume that $\Phi$ preserves positivity.  
 Take $p(x)\in {\mathbb R}[x]$, $p(x)\ge 0$. Then for any $\varepsilon>0$, $p(x)+\varepsilon>0$. Thus $ \Phi(p(x))+\varepsilon\Phi(1)=\Phi (p(x)+\varepsilon )>0$. Taking the limit when $\varepsilon\to 0$ we get that $\Phi(p(x))\ge 0$. 
\end{proof}

The following example shows that,  in general, (3) does not imply (1) and (2).
\begin{ex}
Let $\Phi:{\mathbb R}[x] \to {\mathbb R}[x]$ be defined as follows: 
$\Phi(1)=x^2$, $\Phi(x^i)=0$ for all $i>0$.
Obviously,  $\Phi$ preserves the set of non-negative polynomials but does not preserve the set of positive polynomials since $1$ is mapped to $x^2$ which is only non-negative.
\end{ex}

We are now going to show that in fact this example is in some sense the only possibility, i.e., it essentially describes the whole distinction between positivity- and non-negativity-preservers. 

\begin{theorem}
\label{Lnn1}
Let $\Phi: {\mathbb R}[x]\to {\mathbb R}[x]$  be a  non-negativity-preserver. Then either $\Phi$  is a  positivity-preserver (and therefore an ellipticity-preserver as well)  or $\Phi(1)$ is a polynomial which is only non-negative but not positive. Moreover, in the latter case for any positive polynomial $p(x)\in{\mathbb R}[x]$  the zero locus of $\Phi(p(x))$   is a subset of the zero locus of  $\Phi(1)$.
\end{theorem}

\begin{proof}
Assume that $\Phi$ is a  non-negativity-preserver. Then $\Phi$ sends  positive polynomials to  non-negative ones. Let us assume that $p(x)\in{\mathbb R}[x]$ is positive but its image $\Phi(p(x))$ has real zeros. Since $p(x)$ is positive  there exists $\varepsilon>0$ such that $p(x)-\varepsilon>0$ for all $x$. Thus for its image we have 
$$ \Phi(p(x)-\varepsilon)=\Phi(p(x))-\varepsilon \Phi(1) \ge 0.$$
Set $g(x):=\Phi(p(x))$, $f(x):=\Phi(p(x)-\varepsilon)$, and $h(x):=\Phi(1)$. Since all three polynomials are non-negative and 
$g(x)=\varepsilon h(x)+f(x)$, it follows that for any $x_0$ such that $g(x_0)=0$ one has that $f(x_0)=h(x_0)=0$. 
Since $h(x)$ is a polynomial then either $h(x)\equiv 0$, or $h(x)$ has a finite number of zeros. However, the first possibility is ruled out by Lemma~\ref{L:Phi}, since $h(x)=\Phi(1)$. The second possibility 
implies 
that all positive polynomials whose images are  non-negative but not positive have altogether only a finite number of zeros belonging to the zero locus of $\Phi(1)$. 
\end{proof}

\begin{cor} \label{Lnn}
Let $\Phi: {\mathbb R}[x]\to {\mathbb R}[x]$  be a linear operator such that  $\Phi(1)>0$. Then the conditions (1), (2) and (3) of Theorem~\ref{Lnn2} are equivalent. 
\end{cor}

In exactly the same way  we can show the following. 

\begin{theorem} \label{Lnn_f}
Let $\Phi: {\mathbb R}_n[x]\to {\mathbb R}_n[x]$ be a linear operator  with $\Phi(1)>0$. Then the following conditions are equivalent:
\begin{enumerate}
\item $\Phi$ preserves the set of elliptic polynomials  of degree $\le n$; 
\item either $\Phi$ or $-\Phi$ preserves the set of positive polynomials  of degree $\le n$; 
\item either $\Phi$ or $-\Phi$ preserves the set of non-negative polynomials  of degree $\le n$. 
\end{enumerate}
\end{theorem}

\begin{re} \label{re:3}
Corollary~\ref{Lnn} and Theorem~\ref{Lnn_f} will  allow us to reduce the investigation of non-negativity-, positivity-, and ellipticity-preservers (both in the finite-dimensional and the infinite-dimensional cases) to just one of these three classes of preservers. 
\end{re}


For the sake of completeness notice that for non-linear operators  the situation is different from the one above as the following simple examples show.
\begin{ex} 

1. The bijective map $\Phi_1 : \bR[x] \to \bR[x]$ defined by $$\Phi_1(p(x))=p(x)+c$$ where $c$ is a positive constant, preserves both positivity and non-negativity but does not preserve ellipticity.

2. The bijective map $\Phi_2 : \bR[x] \to \bR[x]$ defined by $$
\Phi_2(p(x))=\begin{cases}  p(x) \qquad \forall p(x)\in \bR[x]\setminus\{x^2+1,-x^2-1\} \\
 -x^2-1 \mbox {  if  }Êp(x)=x^2+1\\ 
x^2+1  \quad\mbox {if  }Êp(x)=-x^2-1
 \end{cases} $$  
 preserves ellipticity, but does not preserve positivity and non-negativity.

3. The bijective map $\Phi_3 : \bR[x] \to \bR[x]$ defined by $$
\Phi_3(p(x))=\begin{cases}Êp(x) \quad\forall p(x)\in \bR[x]\setminus\{\pm x^2\} \\ 
-x^2  \quad \mbox {  if  }Êp(x)=x^2\\  
x^2  \quad\;\;\mbox {  if  }Ê\,p(x)=-x^2 \end{cases} $$
  preserves ellipticity and positivity, but does not preserve non-negativity.
\end{ex}



\section{The  case of diagonal transformations}

The aim of this section is twofold. Firstly we want to recall what was previously known about  positivity- and non-negativity-preservers in the classical case of linear operators acting diagonally in the standard monomial basis of $\bR[x]$ 
and secondly we want to point out  some (known to the specialists in the field, \cite {CC1}, \cite {CC2}) mistakes in the important treatise \cite{I}. 

\subsection{Known correct results}

Let  $T_{\infty}: {\mathbb R}[x]\to {\mathbb R}[x]$ be a linear operator defined by  
\begin{equation}
T(x^i)= \lambda_i x^i \mbox{ for }i=0,1,\ldots
\label{eq:0.1}
\end{equation}
 and, analogously,  let $T_n: {\mathbb R}_n[x]\to {\mathbb R}_n[x]$ be a linear operator  defined by  
\begin{equation}
T_n(x^i)= \lambda_i x^i \mbox{ for } i=0,1,\ldots, n.
\label{eq:0.2}
\end{equation}
 Denote them by $\{\lambda_i\}_{i=0}^{\infty}$ and $\{\lambda_i\}_{i=0}^{n}$, respectively. We will refer to such operators  as {\em diagonal transformations\/}  or {\em diagonal sequences\/}. Diagonal transformations preserving the set of positive polynomials are referred to as {\em $\Lambda$-sequences} in the literature , see \cite{CC1,CC2, K}. 
 Reserving the symbol $\Phi$ for general linear operators we  use  in this section the notation $T\in \Lambda$ to emphasize  that $T$ is a {\em diagonal transformation} preserving positivity.
Multiplying if necessary all elements of our sequence with $-1$, we can assume that $\lambda_0\ge 0$.

\begin{re}
Notice that in the finite-dimensional case we  only need to consider  transformations acting on ${\mathbb R}_n[x]$ for $n$ {\em even}   since there are no positive polynomials of odd degree and a sequence $\{\lambda_i\}_{i=0}^{2k+1}$ preserves the set of positive polynomials in ${\mathbb R}_{2k+1}$ if and only if $\{\lambda_i\}_{i=0}^{2k}$  preserves the set of positive polynomials in ${\mathbb R}_{2k}$.
\end{re}

Let us establish some immediate consequences of the fact that a diagonal transformation $T$ is a positivity-preserver.

\begin{lemma}
Assume that a transformation $T_{\alpha}=\{\lambda_i\}_{i=0}^\alpha$, $\alpha\in {\mathbb N}\cup \{\infty\}$, preserves positivity. Then
\begin{enumerate}
\item $\lambda_i\ge 0$ for any even $i$;
\item $\lambda_i^2\le \lambda_0\lambda_{2i}$ for any $i$.
\end{enumerate}
\end{lemma}

\begin{proof} To settle (1) consider the polynomial $p(x)=x^i+1$ which is positive if  $i$ is even. Thus $T(x)=\lambda_ix^i+\lambda_0$ should be positive as well. Since $\lambda_0>0$, the result follows.

To settle (2) consider the polynomial $p(x)=x^{2i}+ax^i+b$ with  $a^2<4b$. Then $p(x)$ is positive as well as its image $q(x):=\Phi(p(x))=\lambda_{2i} x^{ 2i} + a\lambda_i x^i+b \lambda_0$.  
 If $i$ is odd then the positivity of $q(x)$  is equivalent to the negativity of its discriminant, i.e.,  $D_q:=a^2\lambda_i^2-4b\lambda_0\lambda_{2i}<0$, which implies $\lambda_i^2\le \lambda_0\lambda_{2i}$ since $a^2<4b$. Finally, if $i$ is even then $q(x) $ is positive iff either $D_q<0$ or $D_q>0$ and if  additionally  both roots of the quadratic polynomial 
$\lambda_{2i} z^2 + a\lambda_i z+b \lambda_0$ 
 are negative. In the  first subcase one has  $\lambda_i^2\le \lambda_0\lambda_{2i}$ as for $i$ odd. In the second subcase we obtain that the positive polynomial $x^{2i}-ax^i+b$ is transformed to the  polynomial $\lambda_{2i} x^{2i} - a\lambda_i x^{i}+b \lambda_0$ which has some real roots. To check this notice that  the roots of $\lambda_{2i} z^2 - a\lambda_i z+b \lambda_0$ are opposite to that  of 
 $\lambda_{2i} z^2 + a\lambda_i z+b \lambda_0$ and are, therefore, positive. Thus extracting  their $i$-th root one will get   some positive roots as well. This contradiction shows that the inequality $\lambda_i^2\le \lambda_0\lambda_{2i}$ is necessary  for the positivity-preservation.
\end{proof}

Diagonal transformations which are  positivity-preservers are known to be very closely related to the following class of sequences of real numbers. 

\begin{de}
A sequence $\{\lambda_i\}_{i=0}^{\alpha}$ is called {\em positive definite\/} if for any positive polynomial $p(x)=x^n+a_{n-1}x^{n-1}+\ldots+a_1x+a_0 \in {\mathbb R}_\alpha[x]$ one has that $\lambda_n+a_{n-1}\lambda_{n-1}+\ldots+a_1\lambda_1+a_0\lambda_0>0$, i.e., $T_\alpha(p)(1)>0$.
\end{de}

In the infinite-dimensional case the following characterizations of the set of diagonal   positivity-preservers is known.

\begin{theorem}
\label{LL_PD}  {\rm (\cite{CC1}, \cite[Theorem 1.7]{CC2})}
Let $\{\lambda_k\}_{k=0}^{\infty}$ be a sequence of real numbers. Then the following conditions are equivalent:
\begin{enumerate}\item $\{\lambda_k\}_{k=0}^{\infty}\in \Lambda$;
\item $\{\lambda_k\}_{k=0}^{\infty}$ is a positive definite sequence;
\item $|(\lambda_{i+j})|=\left| \begin{array}{cccc} \lambda_0 & \lambda_1 &\ldots & \lambda_k \\ \lambda_1 & \lambda_2 & \ldots & \lambda_{k+1} \\ \vdots & \vdots & \ddots & \vdots \\ \lambda_k & \lambda_{k+1} & \ldots & \lambda_{2k} \end{array} \right|>0$, for $k=0,1,2\ldots$
\item There exists a non-decreasing function $\mu(t)$ with infinitely many points of increase such that for all $k=0,1,2,\ldots$ one has  
$$
\lambda_k=\int\limits_{-\infty}^{\infty} t^k d\mu(t).  
$$
\end{enumerate}
 \end{theorem}
 
 \begin{proof}   The equivalences 
 (2) $\Leftrightarrow$ (3) $\Leftrightarrow$ (4) are  settled   in \cite[p.~132]{W} independently of condition (1). The implication (1) $\Rightarrow$ (2) is evident so we only have to concentrate on the remaining implication (2) $\Rightarrow$ (1).  Take a positive polynomial  $g(x)=a_lx^l+\ldots + a_1x+a_0$ and set $q(x):=\Lambda(g(x))=\sum\limits_{i=0}^l \lambda_i a_i x^i$.
We want to show that $q(x)> 0$ for all real $x$. 

Condition  (4) implies that $\lambda_i=\int\limits_{-\infty}^{\infty} t^i d\alpha(t)$, where $\alpha(t)$ is a monotone non-decreasing function with infinitely many points of increase. Hence
$$ q(x)=\sum_{i=0}^l \Bigl(\int\limits_{-\infty}^{\infty} t^i d\alpha(t)\Bigr) x^i=\int\limits_{-\infty}^{\infty} \sum_{i=0}^l a_i t^i x^i d\alpha(t)
= \int\limits_{-\infty}^{\infty} g(xt)d\alpha(t) >0 $$
since  $g(xt)>0$ for all $t$. Notice that the above integrals are convergent  for any fixed value of $x$. Thus $q(x)>0$ and the lemma follows. \end{proof}

In the finite-dimensional case one has a similar statement. 
\begin{theorem}
\label{LL_PD_f} \cite{K}
Let $\{\lambda_k\}_{k=0}^{n}$ be a sequence of  real numbers. Then the following conditions are equivalent:
\begin{enumerate}\item $\{\lambda_k\}_{k=0}^{n}\in \Lambda$;
\item $\{\lambda_k\}_{k=0}^{n}$ is a positive definite sequence;
\item $|(\lambda_{i+j})|=\left| \begin{array}{cccc} \lambda_0 & \lambda_1 &\ldots & \lambda_k \\ \lambda_1 & \lambda_2 & \ldots & \lambda_{k+1} \\ \vdots & \vdots & \ddots & \vdots \\ \lambda_k & \lambda_{k+1} & \ldots & \lambda_{2k} \end{array} \right|>0$, for $k=0,1,2\ldots,\frac{n}{2}.$
\item There exists a non-decreasing function $\mu(t)$ with at least $n$ points of increase such that
$$
\lambda_k=\int\limits_{-\infty}^{\infty} t^k d\mu(t)\:, \qquad k=0,1,2,\ldots, n;
$$
\end{enumerate}
 \end{theorem}

\begin{proof}  Repeats  that of Theorem~\ref{LL_PD}.
\end{proof}

\subsection{Known wrong results}

To present some erroneous results from \cite{I}  and the corresponding counterexamples we need  to introduce the following classes of diagonal  transformations.

\begin{de}
We say  that $T_{\alpha}$, $\alpha\in {\mathbb N}\cup{\infty}$, or, equivalently,  the sequence $\{\lambda_i\}_{i=0}^{\alpha}$, is a {\em hyperbolicity-preserver}, if for any hyperbolic $p(x)\in {\mathbb R}_\alpha[x]$ its image $T_\alpha(p(x))$ is hyperbolic. We denote this class of transformations by ${\mathcal H}_\alpha$ or ${\mathcal H}$.
\end{de}

Clearly, this class is the restriction of the earlier defined class of hyperbolicity-preservers to diagonal transformations. 

\medskip
Theorem 4.6.14 of \cite {I} states that  $T\in\Lambda$ if and only if $T^{-1}\in {\mathcal H}$. We will now show that this statement is wrong in both directions.

\begin{st} \label{count}
There exist 

\rm{(i)} $ T\in\Lambda$ such that $T^{-1}\notin {\mathcal H}$;
 
\rm{(ii)} $ T\in{\mathcal H}$ such that $T^{-1}\notin \Lambda$.

\end{st}

\begin{proof}  We present below 3 concrete examples verifying the above claims.
 To illustrate (i) consider  the diagonal transformation $T_4:{\mathbb R}_4[x]\to {\mathbb R}_4[x]$  defined by
the sequence $(\lambda_0,\lambda_1,\lambda_2,\lambda_3,\lambda_4)=\left(\frac {1} {29}, \frac {1}{ 68}, \frac {1}{ 123}, \frac {1}{ 200}, \frac{1}{ 305}\right)$.
By the determinant criteria (3) of Theorem~\ref{LL_PD_f} the operator $T_4$ preserves positivity. However, one can check that its inverse
sends the non-negative polynomial $(x+1)^4$ to the polynomial $
(x+1)(305x^3+495x^2+243x+29)$ possessing two real and two complex roots. \label{ExLnotH}

This example shows  that in the finite-dimensional case there is a diagonal transformation which preserves positivity, but whose inverse does not preserve hyperbolicity. We can extend this example to the infinite-dimensional case as follows. 

 By  \cite[Proposition 3.5]{CC1} that there exists  an infinite sequence $\{\lambda_i\}_{i=0}^\infty\in \Lambda$ such that the sequence of inverses
$\displaystyle\left\{\frac{1}{\lambda_i}  \right\}_{i=1}^\infty\notin {\mathcal H}$. As an explicit example  one can take $\lambda_i=\displaystyle\frac{1}{i^3+5i^2+33i+29}$.

 An  example illustrating (ii) is given in \cite[p. 520]{CC1}, see also \cite[Example 1.8]{CC2}. 
 Namely, the sequence $\{1+i+i^2\}_{i=0}^{\infty}$ corresponds to  a diagonal transformation preserving hyperbolicity. However, the sequence of inverses $\left\{\displaystyle \frac{1}{1+i+i^2}\right\}_{i=0}^{\infty}$ leads to  a diagonal  transformation which is not a positivity-preserver.
\end{proof}

\begin{de}
We say  that a diagonal transformation $T_\alpha$, $\alpha\in {\mathbb N}\cup{\infty}$, generated by the sequence $\{\lambda_i\}_{i=0}^{\alpha}$ is a {\em complex zero decreasing sequence} (CZDS for short), if for any polynomial $p(x)\in {\mathbb R}_\alpha[x]$ the polynomial $T(p)$ has no more non-real  roots (counted with multiplicities) than $p$.  We denote the set of all CZDS by ${\mathcal R}$.  
\end{de}

\begin{re} Obviously, any CZDS preserves hyperbolicity, i.e.,  ${\mathcal R}\subset {\mathcal H}$.  For a while it was  believed  that ${\mathcal R}={\mathcal H}$ until Craven and Csordas found a counterexample \cite {CC1}. Additionally, one can see directly from the  definition that  the inverse of any positive  CZDS is  a $ \Lambda$-sequence, that is, a diagonal positivity-preserver. 
\end{re}

Finally, Theorem 4.6.13 of \cite {I} states that  $T\in\Lambda$ if and only if  $T^{-1}\in {\mathcal R}$, which we disprove below.

\begin{st} 
There exist $T\in\Lambda$ such that $T^{-1}\notin {\mathcal R}$.
\end{st}

\begin{proof} Use the first two counterexamples from  the proof of  Proposition~\ref{count}. 
\end{proof}

\section{Linear ordinary differential operators of finite order}



Our aim in this section is to prove Theorem~\ref{TT1}, i.e., 
to show that there are no positivity-, non-negativity-, and ellipticity-preservers which are linear differential operators of   finite positive order. In fact we are going to show that for any linear differential operator $U$ of order $k\ge 1$ there exists an integer $n$ such that $U:{\mathbb R}_n[x]\to {\mathbb R}_n[x]$ is not a non-negativity preserver. Moreover, we show that one can always choose $n=2k$. Since any positivity-preserver is automatically a non-negativity-preserver and any ellipticity-preserver is a positivity-preserver up to a sign change we will get Theorem~\ref{TT1} in its complete generality from the above statement. 

Denote by $S[s_1,\ldots,s_k]$ the ring of  symmetric polynomials  with real coefficients in the variables $s_1,\ldots, s_k$. Let $\sigma_l$ be the $l$-th elementary symmetric function, i.e., 
 $$\sigma_l=\sum\limits_{j_1<\ldots < j_l} s_{j_1}\cdots s_{j_l}\in S[s_1,\ldots,s_k],\qquad l=1,\ldots,k.$$

\medskip

We will need the following technical fact.

\begin{st} \label{P:w_u}
Let $p(x)=(x-x_1)^2\cdots(x-x_k)^2\in {\mathbb R}[x,x_1,\ldots,x_k]$. Consider the following two families of rational functions:
$$w_l=w_l(x,x_1,\ldots,x_k)=\frac{p^{(l)}(x)}{p(x)},\quad l=1,\ldots,k;$$
\begin{equation} u_l=\sigma_l\left(\frac{1}{x-x_1},\ldots, \frac{1}{x-x_k}\right),\quad l=1,\ldots,k.\label{eq:u-x_i}\end{equation}
Then 

1. $ \displaystyle w_l\in S\left[\frac{1}{x-x_1},\ldots, \frac{1}{x-x_k}\right]$, $l=1,\ldots,k$. 

\medskip

2. For any $l=1,\ldots, k $ one has that
 $$w_l=2l! u_l + g_l(u_1,\ldots,u_{l-1}),$$
  where  $g_l\in {\mathbb R}[y_1,\ldots,y_{l-1}]$, $l=1,2,\ldots,k$, are certain polynomials (that can be found explicitly but we will not need their explicit form;  in particular,  $g_1\equiv 0$).
\end{st}


\begin{proof} 
%
%
%
Set $p_1(x)=(x-x_1)\cdots(x-x_k)$. Then one can immediately check that 
$$\frac{p_1^{(i)}(x)}{p_1(x)}=i!u_i$$
for all $i=1,\ldots,k$. Using the Leibniz rule we get 
$$
p^{(l)}(x)=(p_1^2(x))^{(l)}=2p_1^{(l)}(x)p_1(x)+\sum\limits_{i,j\ge 1, i+j=l} c_{i,j} p_1^{(i)}(x)p_1^{(j)}(x), 
$$
 where $c_{i,j}\ge 0$ are certain binomial coefficients.
Thus
$$
w_l=\frac{p^{(l)}(x)}{p(x)}= 2\frac{p_1^{(l)}(x)}{p_1(x)}+\sum\limits_{i,j\ge 1, i+j=l} c_{i,j} \frac{p_1^{(i)}(x)}{p_1(x)} \cdot \frac{p_1^{(j)}(x)}{p_1(x)}  = 2l!u_l+\sum\limits_{i,j\ge 1, i+j=l}c_{i,j}u_i u_j.$$
 The result follows. \end{proof}

We are now ready to prove the first main result of this paper.  
\begin{theorem} \label{T1}
Let $U_Q:{\mathbb R}[x]\to {\mathbb R}[x]$ be a linear ordinary differential operator of order $k\ge 1$ of the form  (\ref{BasOper}) with polynomial coefficients $Q=(q_0(x), q_1(x),\ldots, q_k(x))$,   $q_i\in {\mathbb R}[x]$, $q_k(x)\not\equiv 0$. 
Then for any such coefficient sequence $Q$  the operator $U_Q$ does not preserve the set of non-negative polynomials of degree~$2k$.
\end{theorem}

\begin{proof} We assume that $U_Q \not \equiv 0$. Since $U_Q(1)=q_0(x)$ an obvious necessary condition for the operator $U_Q$ to preserve non-negativity is that $q_0(x)$ itself is a non-negative polynomial. Moreover,  $q_0(x)$ does not vanish identically by Lemma~\ref{L:Phi}. 

We will now construct a  non-negative polynomial 
$$
p(x):=p_{x_1,\ldots,x_k}(x)= (x-x_1)^2(x-x_2)^2\cdots (x-x_k)^2 \in {\mathbb R}[x,x_1,\ldots,x_k].
$$
such that its image under the action of $U_Q$ attains negative values. 
For this we define
$$
R(x,x_1,\ldots,x_k)=\frac {U_Q (p(x))}{p(x)} = 
q_0(x)+q_1(x)\frac{p'(x)}{p(x)}+q_2(x)\frac{p''(x)}{p(x)}+\ldots +q_k(x)\frac{p^{(k)}(x)}{p(x)}.
$$

Then in the notation of Proposition~\ref{P:w_u} we have that
$ R(x,x_1,\ldots,x_k)=q_0(x)+q_1(x)w_1+\ldots+q_k(x)w_k$. Let us fix $x_0\in {\mathbb R}$ such that $x_0\ne 0$ and for any $i=0,1,\ldots,k$ either $q_i(x)\equiv 0$ or $q_i(x_0)\ne 0$. Set $\alpha_i=q_i(x_0)$, $i=1,\ldots,k$, and $r(x_1,\ldots, x_k)=R(x_0,x_1,\ldots,x_k)$. Then $r(x_1,\ldots, x_k)$ is a linear form in $w'_1,\ldots,w'_k$, where $w_i'=w_i(x_0)$. Thus there exist $a_1,\ldots,a_k\in {\mathbb R}$ such that $q_0(x_0)+
\sum\limits_{i=1}^k\alpha_ia_i<0$. (Notice that by our choice of $x_0$ one has  $q_0(x_0)>0$.) 


Let now $b_1,\ldots,b_k\in {\mathbb R}$ be defined by
\begin{equation}
b_1=a_1,\quad b_i=\frac{1}{2i!}(a_i-g_i(b_1,\ldots,b_{i-1})),\quad i=2,\ldots,k,
\label{eq:a_b}
\end{equation}
where $g_i$ are defined in Proposition~\ref{P:w_u}.
Consider  the system of equations 
\begin{equation}
b_i =\sigma_i\left(\frac{1}{x_0-t_1},\ldots, \frac{1}{x_0-t_k}\right), \quad i=1,\ldots,k,
\label{eq:b_t}
\end{equation}
with unknowns $t_1,\ldots,t_k$. It follows from the Vi\`eta theorem that the $k$-tuple $(t_1,\ldots,t_k)$ solves  system (\ref{eq:b_t}) if and only if $\frac{1}{x_0-t_1},\ldots, \frac{1}{x_0-t_k}$ are the roots of the equation 
$$
z^k-b_1z^{k-1}+b_2z^{k-2}-\ldots \pm b_{k-1}z\mp b_k=0.
$$
Note that since $b_i\in {\mathbb R}$ the roots of the latter equation are either real or complex conjugate. Wlog  we can always assume that they are ordered so that  $z_1=\overline{z_2},\ldots, z_{2i-1}=\overline {z_{2i}}\in {\mathbb C}$, $z_{2i+1},\ldots,z_k\in {\mathbb R}$. Thus the $k$-tuple $(t_1,\ldots,t_k)$, where $t_i=x_0-\frac{1}{z_i}$, $i=1,\ldots, k$, solves  system  (\ref{eq:b_t}). Obviously,  $t_1=\overline{t_2},\ldots, t_{2i-1}=\overline {t_{2i}}\in {\mathbb C}$, $t_{2i+1},\ldots,t_k\in {\mathbb R}$. 


Let us substitute $x=x_0$, $x_i=t_i$, $i=1,\ldots,k$, in the functions $w_l$, $l=1,\ldots,k$, defined in Proposition \ref{P:w_u}. 
By the definition of the $u_i$'s and using the fact that the $t_i$'s solve system~(\ref{eq:b_t}) we get that $u_i(x_0,t_1,\ldots,t_k)=b_i$, $i=1,\ldots,k$. Thus (\ref{eq:a_b}) implies that  
$$w_i(x_0,t_1,\ldots,t_k)=a_i,\quad  i=1,\ldots,k.$$
 Hence $R(x_0,t_1,\ldots,t_k)<0$ by the choice of the $a_i$'s. Thus $U_\alpha(p_{t_1,\ldots,t_k})(x_0)<0$ since $p_{x_1,\ldots,x_k}(x)\ge 0$ for all $x\in {\mathbb R}$ and all $x_i, i=1,\ldots,k$. This  contradiction proves Theorem~\ref{T1}.
\end{proof}

As we mentioned above Theorem~\ref{T1} together with the results of \S 2 imply Theorem~\ref{TT1}. 

\medskip 

We are now going to strengthen Theorem~\ref{T1} and show that   wide subclasses of  linear ordinary differential operators of finite order $k$ do not  preserve non-negativity even in degrees much smaller than $2k$. In particular,  the next statement shows that no linear differential operator of odd order $k$ preserves the set of non-negative polynomials in ${\mathbb R}_{k+1}[x]$. 

\begin{st} \label{T2}
Let $U_Q:{\mathbb R}[x]\to {\mathbb R}[x]$ be a linear differential operator of odd order $k\ge 1$ of the form  (\ref{BasOper})  with polynomial coefficients $Q=(q_0(x),q_1(x),\ldots, q_k(x))$, $q_i(x)\in {\mathbb R}[x]$, $i=1,\ldots,k$, $q_k(x)\not\equiv 0$. 
Then $U_Q$ does not preserve the set of non-negative polynomials of degree smaller than or equal to $k+1$.
\end{st}

\begin{proof} 
Consider the polynomial $p_t(x)=(x-t)^{k+1}$.  It is non-negative since $k$ is odd. Note that
$$ r_Q(x,t) := \frac{U_Q(p_t(x))}{p_t(x)}= q_0(x) +(k+1) q_1(x) \frac{1}{x-t} +\ldots + (k+1)! q_k(x)\frac{1}{(x-t)^k}$$
and set $u:=\frac{1}{x-t}$. Then the function 
$$ g(x,u):= q_0(x) +(k+1) q_1(x)u+\ldots + (k+1)!q_k(x)u^k$$ 
is a polynomial in $u$. Fixing $x_0$ such that $q_k(x_0)\ne 0$ we obtain that $g(x_0,u)$ is a polynomial in $u$ of odd degree. Hence there exists $u_0$ such that $g(x_0,u_0)<0$. Now for $t_0=x_0-\frac{1}{u_0}$ we get that $r_Q(x_0,t_0)<0$. Thus $U_Q(p_{t_0})(x_0)<0$ since $p_t(x)\ge 0$ for all $(x,t)$. This  contradiction finishes the proof. 
\end{proof}

The next result shows that there is also a  large class of linear differential operators of even order $k$ which does not preserve non-negativity in ${\mathbb R}_k[x]$.

\begin{st} \label{T3}
Let $U_Q:{\mathbb R}[x]\to {\mathbb R}[x]$ be a linear differential operator of even order $k$ of the form  (\ref{BasOper}) with polynomial coefficients $Q=(q_0(x),q_1(x),\ldots, q_k(x))$, $q_i(x)\in {\mathbb R}[x]$, $i=1,\ldots,k$, $q_k(x)\not\equiv 0$. 
Assume in addition that either there exists $x_0\in{\mathbb R}$ such that $q_k(x_0) < 0$ or  there exists $x_0\in{\mathbb R}$ such that  $q_k(x_0)=0$ and $q_{k-1}(x_0)\ne 0$. Then $U_Q$ does not preserve the set of non-negative polynomials of degree smaller  than or equal to~$k$.
\end{st}

\begin{proof} 
The polynomial $p_t(x)=(x-t)^{k}$  
is non-negative since $k$ is even. 
Similar to the above one has 
$$ r_Q(x,t) := \frac{U_Q(p_t(x))}{p_t(x)}= q_0(x) +k q_1(x) \frac{1}{x-t} +\ldots + k! q_k(x)\frac{1}{(x-t)^k}.$$

As before we set $u:=\frac{1}{x-t}$ and consider the function 
$$ g(x,u):= q_0(x) +k q_1(x)u+\ldots + k! q_k(x)u^k,$$ 
which is a polynomial in $u$.
If there exists $x_0$ such that $q_k(x_0)<0$ then $g(x_0,u)$ is a polynomial in $u$ which is negative for sufficiently large values of $u$. 
If $q_k(x_0)=0$ then $g(x_0,u)$ is a polynomial in $u$ of odd degree. In both cases there exists $u_0$ such that $g(x_0,u_0)<0$. Now for $t_0=x_0-\frac{1}{u_0}$ we get that $r_Q(x_0,t_0)<0$. Thus $U_Q(p_{t_0})(x_0)<0$ since $p_t(x)\ge 0$ for all $(x,t)$. This  contradiction accomplishes the proof. 
\end{proof}


\begin{cor} \label{CT3}
Let $U_Q:{\mathbb R}[x]\to {\mathbb R}[x]$ be a linear differential operator of order $k$ of the form  (\ref{BasOper}) with polynomial coefficients $Q=(q_0(x),q_1(x),\ldots, q_k(x))$, $q_i(x)\in {\mathbb R}[x]$, $i=1,\ldots,k$, $q_k(x)\not\equiv 0$. 
Assume that there exists an even integer  $i\in \{2,4,\ldots,k\}$ such that either there exists $x_0\in{\mathbb R}$ such that  $q_i(x_0)< 0$ or  there exists $x_0\in{\mathbb R}$ such that $q_i(x_0)=0$ and $q_{i-1}(x_0)\ne 0$.
Then $U_Q$ does not preserve non-negativity in ${\mathbb R}_l[x]$ for any $l\ge i$.
\end{cor}

\begin{proof}  If $l\ge i$ and $U_Q$ preserves non-negativity in ${\mathbb R}_l[x]$ then $U_Q$ preserves non-negativity in ${\mathbb R}_i[x]$. The restriction of $U_Q$ to the space ${{\mathbb R}_i[x]}$ is given by 
$$ q_0(x)+q_1(x) \frac{d}{dx}+q_2(x)\frac{d^2}{dx^2}+\ldots + q_i(x) \frac{d^i}{dx^i}$$ and the result follows from Theorem~\ref{T3}. 
\end{proof}

\begin{re} In particular, if $i=2$ in Corollary~\ref{CT3} then there is no degree $l$, $l\ne 0$, such that $U_Q:{\mathbb R}_l[x]\to {\mathbb R}_l[x]$ preserves non-negativity. Thus for a ``generic''  linear differential operator $U_Q: {\mathbb R}[x]\to {\mathbb R}[x]$ with non-constant coefficients there is no $l$ such that $U_Q$ preserves non-negativity in ${\mathbb R}_l[x]$.
 However, the following example shows that there are  linear differential operators with non-constant coefficients which preserve positivity on~$\bR_k[x]$ for any even~$k$.
\end{re} 

\begin{ex} \label{Ex:exists}
Let $k$ be an even positive integer and consider the linear differential operator of order $k$ given by $$U_k^{a,b}=1+(ax^k+b)\frac{d^k}{dx^k}:{\mathbb R}_k[x]\to {\mathbb R}_k[x].$$  Then $U_k^{a,b}$ preserves positivity on $\bR_k[x]$ if $a\ge 0$, $b\ge 0$.
\end{ex}
Indeed, $U_k^{a,b}(p(x))=p(x)+k!a_k(ax^k+b)\ge p(x)>0$ for any positive polynomial $p(x)=a_kx^k+\ldots \in \bR_k[x]$.

\section{Linear ordinary differential operators with constant coefficients}

In this  section we will prove  Theorem~\ref{TC:Fd}. Take a sequence $\alpha=(\alpha_0,\alpha_1,\ldots, \alpha_k)$  of real numbers. Denote by $U_\alpha$ the following  linear differential operator of order $k$ 
\begin{equation}\label{CCF}
U_\alpha=\alpha_0+\alpha_1\frac{d}{dx}+\ldots +\alpha_{k}\frac{d^{k}}{dx^{k}}
\end{equation}
 with constant coefficients.

By Theorem~\ref{T1} there are no finite order linear differential operators on ${\mathbb R}[x]$ preserving  positivity. 
However,  in the case of polynomials of bounded degree, i.e., belonging to  the finite-dimensional space $\bR_k[x]$, there are such linear differential operators, see Example~\ref{Ex:exists}.


%


Theorem~\ref{TC:Fd} follows easily from the next statement of Remak \cite{Re} and Hurwitz \cite{Hu} which for the sake of completeness  we present with its  proof.

\begin{theor} \label{ReHu} For an even integer $k=2l$  and a     sequence of real numbers $\alpha=(\alpha_0,\alpha_1,\ldots,\alpha_k )$ consider  the   linear ordinary differential operator  (\ref{CCF}) 
with constant coefficients. Then the operator $U_{\alpha}$ 
 preserves positivity (resp., non-negativity) in $\bR_k[x]$ if and only if one of the following two equivalent conditions holds:
 \begin{enumerate}
 \item   for any positive (resp., non-negative) polynomial $p(x)=a_{k}x^{k}+\ldots + a_1x+a_0$ one has  that 
 $$ U_\alpha(p(0))=a_0\alpha_0+a_1\alpha_1+\ldots + k! a_{k} \alpha_{k}>0 \ (\mbox{resp., }\ge 0); 
 $$
 \item 	the following $(l+1)\times(l+1)$  Hankel matrix 
 $$ \begin{pmatrix}  \alpha_0 & 1!\alpha_1 & 2!\alpha_2 & \ldots & l!\alpha_l \\ 1!\alpha_1 & 2!\alpha_2 & 3!\alpha_3 & \ldots & (l+1)!\alpha_{l+1} \\  2!\alpha_2 & 3!\alpha_3& 4!\alpha_4 & \ldots & (l+2)!\alpha_{l+2}  \\ \vdots & \vdots & \vdots & \ddots& \vdots \\ l!\alpha_l & (l+1)!\alpha_{l+1} & (l+2)!\alpha_{l+2} & \ldots & (2l)!\alpha_{2l}  \end{pmatrix}$$ 
is positive definite (resp., positive semi-definite).    
 \end{enumerate}
 \end{theor}

We start with  the following observation. 

\begin{lemma} \label{L:shift}
The operator  $U_{\alpha}:{\mathbb R}_{k}[x]\to {\mathbb R}_{k}[x]$   of the form (\ref{CCF}) commutes with shifts of the independent  variable $x$. In other words,  for any polynomial $p(x)\in {\mathbb R}_k[x]$ set $q(x)=U_\alpha(p(x))$. Then for any $x_0\in {\mathbb R}$ we have that $q(x-x_0)=U_{\alpha}(p(x-x_0))$.
\end{lemma}

\begin{proof}
Take any $p(x)=a_lx^l+\ldots+a_0,\; l\le k$. Then for any positive integer $i$ we have that $p^{(i)}(x)=a_l\frac{l!}{i!}x^{l-i}+\ldots + a_ii!$. Thus $(p^{(i)})(x-x_0)=(p(x-x_0))^{(i)}$. Since the coefficients of $U_{\alpha}$ are constant, the result follows.
\end{proof}

 

 
 \begin{proof}[Proof of Theorem~\ref{ReHu}]   The equivalence between the conditions  (1) and (2) in the formulation of  Theorem~\ref{ReHu} is exactly the same fact as  the equivalence between (2) and (3) in Theorems \ref{LL_PD}-\ref{LL_PD_f} for $\lambda_0=\alpha_0$, $\lambda_i=i!\alpha_i$, $i=1,\ldots,k$. It is proven in the required generality in Theorem 11, 
 p. 133 of \cite{W}.  What we need is to show that the assumption that $U_\alpha$ is a non-negativity-preserver (resp., a positivity-preserver) in $\bR_k[x]$  is equivalent to  condition (1). Indeed,  if $U_\alpha$ preserves non-negativity in $\bR_k[x]$, then $a_0\alpha_0+a_1\alpha_1+\ldots + k! a_{k} \alpha_{k}=U_{\alpha} (p(0))\ge 0 $ for any non-negative polynomial $p(x)\in \bR_k[x]$. Assume now that for any non-negative polynomial $p(x)$ one has that $a_0\alpha_0+a_1\alpha_1+\ldots + k! a_{k} \alpha_{k}\ge 0 $. Set $q(x):=U_Q(p(x))$. By assumption we have that $q(0)\ge 0$ and we want to show that $q(x)$ is non-negative. For any $x_0\in {\mathbb R}$  consider $g_{x_0}(x):=q(x+x_0)$. By Lemma~\ref{L:shift} we have that
 $g_{x_0}(x)=U_\alpha (p(x+x_0))$, but $f(x):=p(x+x_0)$ is a non-negative polynomial. Thus by condition (1) one has $U_\alpha(f(0))\ge 0$, i.e., 
 $$q(x_0)=g_{x_0}(0)=U_\alpha(f(0))\ge 0$$ 
 for any $x_0\in {\mathbb R}$. Simple additional consideration shows that the same argument with the strict inequality in condition (1) works for positivity-preservers. 
\end{proof}

\begin{re}
 Theorem~\ref{ReHu} provides the classification of linear differential operators with constant coefficients of an even order $k$ which preserve positivity in ${\mathbb R}_k[x]$. On the other hand, by Theorems~\ref{TT1} and \ref{T1} there are no linear differential operators with constant coefficients of even order $k$  that preserve positivity in ${\mathbb R}_{2k}[x]$. Below we bridge this gap between $k$ and $2k$ for operators with constant coefficients by showing that there are no such operators  of order $k$  that preserve positivity (or non-negativity, or ellipticity) in ${\mathbb R}_{l}[x]$ for any~$l > k$.
 \end{re}

\begin{st}
Let $k$ be a positive integer and let  $\alpha=(\alpha_0, \alpha_1,\ldots,\alpha_k)$ be a sequence of real numbers. Consider the operator $U_{\alpha}$ of the form (\ref{CCF}). 
 Then for any $l > k$ the operator $U_\alpha: \bR_l[x]\to \bR_l[x]$ does not preserve positivity. 
\end{st}

\begin{proof}
Wlog   we can assume that $l$ is even. We can also assume $\alpha_0>0$ and  at least one more entry $\alpha_j$ in the sequence $(\alpha_0, \alpha_1,\ldots,\alpha_k)$ is non-vanishing. (The cases when either $\alpha_0\le 0$ or only $\alpha_0$ is non-vanishing are trivial.) Take any (not necessarily positive!) polynomial $p(x)=a_kx^k+...+a_1x+a_0$ of degree at most $k$ such that $a_0>0$ and $U_\alpha(p(0))=a_0\alpha_0+a_1\alpha_1+...+k!a_k\alpha_k<0$.  Since both $\alpha_0$ and $\alpha_j$ are non-vanishing such a $p(x)$ always exists. Consider now $P(x)=Mx^l+p(x)$ where $M$ is a large positive constant. By our assumptions one can always choose such a large $M$ that $P(x)$ becomes 
positive. At the same time $U_\alpha(P(0))=U_\alpha(p(0))<0$. The latter contradicts to the condition (1) of Theorem~\ref{ReHu} implying that $U_\alpha$ does not preserve positivity in $ \bR_l[x]$. 
\end{proof}

\begin{re} Finally, notice that in order to get Theorem~\ref{TC:Fd} from Theorem~\ref{ReHu} it suffices to observe that for any positive even integer $k$ the action of the operator $U_\alpha$  of infinite order of the form (\ref{ConsOper}) on the space $\bR_k[x]$ coincides with the action of its truncation (\ref{CCF}). 
\end{re}


\begin{thebibliography}{999}
\bibitem{BB}
J.~Borcea, P.~Br\"and\'en, {\em Lee-Yang theory and linear operators preserving
stability}, preprint (2008). 

\bibitem{BBS1} 
J.~Borcea, P.~Br\"and\'en, B.~Shapiro, {\em Classification of hyperbolicity and stability preservers: the multivariate  Weyl algebra case},  arXiv:math.CA/0606360.

\bibitem{BBS2} 
J.~Borcea, P.~Br\"and\'en, B.~Shapiro, {\em P\'olya-Schur master theorems for circular domains and their boundaries}, to appear in Annals of Mathematics,  arXiv:math.CV/0607416. 

\bibitem{CC1} 
T.~Craven, G.~Csordas, {\em Problems and theorems in the theory of multiplier sequences}, Serdica Math. J. {\bf 22} (1996), 515--524.

\bibitem{CC2} 
T.~Craven, G.~Csordas, {\em Complex zero decreasing sequences}, Methods Appl. Anal. {\bf 2} (1995), 420--441.

\bibitem{Ehr}
R.~Ehrenborg, {\em  The Hankel determinant of exponential polynomials}. Amer. Math. Monthly {\bf 107}  (2000), no. 6, 557--560.

\bibitem{Hu} 
A.~Hurwitz, {\em \"Uber definite Polynome}, Math. Ann. {\bf 73} (1913), 173--176.


\bibitem{I} 
L.~Iliev, Laguerre entire functions. 2nd ed., Publ. House of the Bulgarian Acad. Sci., Sofia, 1987, 188 pp.

\bibitem{K} 
M.~D.~Kostova, {\em \"Uber die $\lambda $-Folgen} (German) [On $\lambda $-sequences], C. R. Acad. Bulgare Sci. {\bf 36} (1983), 23--25. 
 
\bibitem{Lom} 
J.~S.~Lomont, J.~Brillhart, Elliptic polynomials. Chapman \& Hall/CRC, Boca Raton, FL, 2001, xxiv+289 pp.

\bibitem{Mo} 
T.~S.~Motzkin, {\em Algebraic inequalities}, in ``Proc. Sympos. Wright-Patterson Air Force Base,'' Ohio, 1965, pp.~199--203, Academic Press, New York.

\bibitem{Pak} 
F.~B.~Pakovich, {\em Elliptic polynomials} (Russian), Uspekhi Mat. Nauk. 
{\bf 50} (1995), 203--204; English translation in Russian Math. Surv. {\bf 50} (1995), 1292--1294.

\bibitem{PD} 
A.~Prestel, C.~N.~Delzell, Positive polynomials. From Hilbert's 17th problem to real algebra. Springer Monographs in Mathematics, Springer-Verlag, Berlin, 2001, viii+267 pp.

\bibitem{PS}
G.~P\'olya, I.~Schur,  {\em \"Uber zwei Arten von
Faktorenfolgen in der Theorie der algebraischen Gleichungen}, J. Reine
Angew. Math. {\bf 144} (1914), 89--113.

\bibitem{PSz} 
G.~P\'olya, G.~Szeg\"o, Problems and theorems in analysis. II. Theory of functions, zeros, polynomials, determinants, number theory, geometry. Translated from the German by C. E. Billigheimer. Reprint of the 1976 English translation. Classics in Mathematics. Springer-Verlag, Berlin, 1998. xii+392 pp.

\bibitem{Re} 
R.~Remak, {\em Bemerkung zu Herrn Stridsbergs Beweis des Waringschen Theorems}, 
Math. Ann. {\bf 72} (1912), 153--156. 

\bibitem{W} 
D.~V.~Widder, The Laplace Transform. Princeton Math. Series Vol 6, Princeton Univ. Press, Princeton, NJ, 1941, x+406 pp.

\end{thebibliography}
\end{document}